\documentclass[12pt]{article}
\usepackage{amsmath}
\usepackage{amssymb}
\usepackage{graphicx}
\newtheorem{theorem}{Theorem}

\newtheorem{lemma}{Lemma}

\newtheorem{definition}{Definition}
\newtheorem{proposition}{Proposition}

\newenvironment{proof}[1][Proof]{\textbf{#1.} }{\ \rule{0.5em}{0.5em}}
\everymath{\displaystyle}

\begin{document}

\centerline{\textbf {New bounds for the b-chromatic number of vertex deleted
    graphs}}

\vspace{8.0mm}
\centerline{\textbf{Renata DEL-VECCHIO and Mekkia KOUIDER}}
%\centerline{\small {emails: renata@vm.uff.br, \hspace{1.0mm} km@lri.fr}}

\vspace{15.0mm}
\noindent \textbf{Abstract.} A b-coloring of a graph is a proper coloring of its vertices such
that each color class contains a vertex adjacent to
at least one vertex of every other color class. The b-chromatic number of a
graph is the largest integer $k$ such that the graph
has a b-coloring with $k$ colors. In this work we present lower bounds for
the b-chromatic number of a vertex-deleted subgraph of a graph,
particularly regarding two important classes, quasi-line and chordal graphs.We
also get bounds for the b-chromatic number of $G-\{x\}$, when $G$ is a graph
with large girth.
\\

\vspace{1.0mm}


\vspace{1.0mm}


\vspace{2.0mm}

\section{Introduction}

\vspace{2.0mm}

    All graphs considered in this work are simple and undirected.
Let $ G= (V(G),E(G))$ be an undirected graph where $V(G)$ and $E(G)$ are the sets of its vertices and
edges, respectively. If
$A \subset V(G)$, we denote by $\langle A\rangle$ the induced subgraph
generated by $A$. For any vertex $x$ of a graph G,
the neighborhood of $x$ is the set $N(x) = \{y \in V(G) | xy \in E(G) \}$. The degree of a vertex $x$ is the cardinality of $N(x)$ and it is denoted by $d(x)$.  We denote by $\Delta(G)$ the maximum degree of $G$. If $x,y \in V(G)$, the distance between  $x$ and $y$ (that is, the length of the shortest $x-y$-path) is represented by $d(x,y)$ and $N^2(x)= \{y \in V(G) | d(x,y)=2\}$.
The {\it girth} of $G$ is the length of its shortest cycle, and is
denoted by $g(G)$. 

Let $G$ be a graph with a proper vertex coloring. Let us denote by $C_i$ the set of vertices of color $i$, herein called the class of color $i$. Let $x_i$ denote a vertex $x$ of color $i$; $x_i$ is said a \textit{color-dominating vertex} (or,
\textit{b-dominating vertex}) if $x_i$ is adjacent to at least one vertex in each of the other classes. A color $i$ is a \textit{dominating color} if there is at least one vertex $x_i$ that is color-dominating.
If $y_i$ is a vertex which is not color-dominating, at least one color $j \neq i$
does not appear in $N(y_i)$. The color $j$ is said a \textit{missing color} in
$N(y_i)$ or simply a missing color of $y_i$.

A {\it b-coloring} is a proper coloring of its vertices
such that each color class contains a color-dominating vertex. 

The b-chromatic number $b(G)$ is the largest integer  $k$ such that $G$ admits a b-coloring with $k$ colors.
Since this parameter has been introduced by R. W. Irving
and D. F. Manlove \cite{Irving}, it aroused the interest of many researchers
as we can see in \cite{Hoang}, \cite{Bonomo}
and \cite{Ana} and, more recently, in \cite{Jak}.

For a vertex $x$ of $G$, let $G-\{x\}$ be the vertex-deleted subgraph of $G$
obtained by deleting $x$ and all edges incident to $x$.
It is known that the chromatic number of $G-\{x\}$ can have a maximum
variation of one unit compared to the chromatic number of $G$.
However, this is not true for the b-chromatic number - the difference between
$b(G)$ and $b(G-\{x\})$ can be arbitrarily large.
This fact motivates the search for bounds to the b-chromatic number
(see \cite{Bala} and \cite{Kouider}).
For general graphs S.F.Raj and R.Balakrishnan proved that:

\begin{theorem}{\cite{Bala}}
For any connected graph of order $n\geq 5$, and for any vertex $x \in V(G)$,

$b(G)-(\lceil n/2 \rceil -2) \leq b(G-\lbrace x\rbrace) \leq b(G)+(\lfloor n/2 \rfloor -2)$
\end{theorem}
The bounds are sharp.

\vspace{2.0mm}

In \cite{KouiderM}, some upper bounds of $b(G -\lbrace x\rbrace)$ have been established in
some classes of graphs as quasi-line graphs, graphs of large girth and  chordal graphs. A {\it chordal} graph is a graph such that every cycle of length at least $4$ has a chord.
A graph is a {\it quasi-line} graph if the neighborhood of each vertex is covered by at most two cliques. In particular, claw-free graphs (i.e. graphs without induced $K_{1,3}$) are quasi-line graphs.
The following results have been shown.

\begin{theorem}\cite{KouiderM} \label{quasi}
  Let $G=(V,E)$ a graph.

 1) If $G$ is a quasi-line-graph, then for each vertex $x$, $$b(G-\lbrace x\rbrace) \leq b(G)+2$$

 2) If $G$ is any graph of girth at least $5$, then for each vertex $x$, $$b(G-\lbrace x\rbrace) \leq b(G)+1$$

\end{theorem}
\vspace{4mm}

\begin{theorem} \cite{KouiderM} \label{coch}
Let $G=(V,E)$ be a chordal graph of clique-number $\omega$
and b-chromatic number $b(G)$.\
Then, for each vertex $x$,
$$b(G-\lbrace x\rbrace) \leq  b(G)+1 + \sqrt{d(x)-1}$$
$$b(G-\lbrace x\rbrace) \leq  b(G)+1 + \sqrt{\omega -1}$$

\end{theorem}
In this work we present lower bounds for $b(G-\{x\})$ in terms of $b(G)$,
particularly regarding two important classes of graphs, quasi-line and chordal
graphs. We also obtain a lower bound for $b(G-\{x\})$ for graphs of large girth.

Besides this introduction we have three more sections. In the second
one we present a lower bound for $b(G-\{x\})$, when $G$ is a general graph
and another bound for quasi-line graphs.
The third section is devoted to the study of chordal graphs, obtaining also
a lower bound for $b(G-\{x\})$ in this class. Finally, in the last section,
we analyse graphs with girth at least $5$, presenting also here a lower bound
for $b(G-\{x\})$.

\section{General bound and quasi-line graphs}
We begin this section with a lower bound for the b-chromatic number of a
vertex deleted subgraph of any graph.

\begin{proposition}
  For every vertex $x \in V(G)$, $b(G-\{x\})\geq b(G)-d(x).$
  
\end{proposition}

\proof
Let $x \in V(G)$ be a fixed vertex. For a b-coloring of $G$, let $i$ be the color of $x$. We consider
two cases. First suppose 
that  each color has at least one color-dominating vertex in $G-\{x\}$; then the b-coloring of $G$ is also
a b-coloring of $G-\{x\}$, so $b(G-\{x\})\geq b(G)$.
Now, let us consider that there is a color with no color-dominating vertices
in  $G-\{x\}$. We have then two possibilities:
\begin{itemize}
\item There is no color-dominating vertex of color $i$ in $G-\{x\}$, that is, $C_i$ has no color-dominating vertex then; for each vertex $z$ in $C_i$, there is at least one color missing in $N(z)$. We can change
  the color of each vertex $z$ in $C_i$ by a missing color in $N(z)$, eliminating the color $i$.As $C_i$ is a stable set, the new coloring is proper.
  For this case we have $b(G-\{x\})\geq b(G)-1 \geq b(G)-d$.
  
%The precedent case is now excluded.
\item There is a vertex $y \in N(x)$ such that $y$ was the color-dominating vertex of color $s, s \neq i$ in $G$ and there is no more color-dominating vertex
of color $s$ in $G-\{x\}$. As $C_s$ has no color-dominating vertices, we then change the color $s$ of $y$ by $i$ and, for each other vertex $z$  in $C_s$, we change the color $s$ for a missing color of $z$, eliminating the color $s$. As $C_s$ is a stable set, the new coloring is proper. We repeat this process for
all vertices in $N(x)$ in the same conditions as $y$.
We do this for at most $d(x)$ vertices. 
  In this case we obtain $b(G-\{x\})\geq b(G)-d(x)$.
\end{itemize}

\begin{flushright}
$\blacksquare$
\end{flushright}

If $\Delta(G) < \left\lceil\frac{n}{2}\right\rceil -2 $, this bound is better than the lower bound in \cite{Bala}.

Note that there exist chordal (resp. quasi-line) graphs $G$ such that
$b(G-\{x\})$ is stricly less than $b(G)$.
For example, let $G_0$ be a chordal graph obtained from a chordal graph $H$ and a new vertex $x$ joined to every vertex of $H$. Then
$b(G_0~-x)= b(G_0)-1$.\\

\begin{theorem}
  If $G$ is a quasi-line graph then, for every vertex $x \in V(G)$,
$b(G-\{x\})\geq b(G)-2$.
\end{theorem}

\proof
Let $x \in V(G)$ be a fixed vertex. 
%As $G$ is a quasi-line graph, there are no three mutually non %adjacent vertices in $N(x)$. 
$N(x)$ is covered by at most two cliques $K_1$ and $K_2$. Let $i$ be the color of $x$.

Considering a b-coloring of $G$, there is at most two vertices
$u_k, u'_k \in N(x)$ with the same color $k$. Again, by the
fact that $G$ is quasi-line, $N^{2}(x) \cap N(u_k)$ is a clique as it is
independent from the neighbour $x$ of $u_k$. Analogously
$N^2(x) \cap N(u'_k)$ is a clique.\\
We delete the vertex $x$.If in $G-\{x\}$ each color is dominating, then
$b(G-\{x\}) \geq b(G).$
Let $i$ be the color of $x$in $G$ . We may suppose that $G-\{x\}$ has a color-dominating vertex of color $i$ otherwise we color each
vertex of color $i$ by a missing color and we get a b-coloring of $G-\{x\}$ by
$b(G)-1$ colors .\\

We choose a color $s$ that is no more dominating, which means that there is no more color-dominating vertices of color $s$. Each vertex of color $s$ has a
missing color.

If the color $s$ had more than one color-dominating vertex in $G$, then it had
exactly two color-dominating vertices $w_s \in K_1$ and $w'_s \in K_2$. We recolor both of them by $i$. We then recolor each other vertex of color $s$ by a missing color. In this way we obtain a b-coloration of $G-\{x\}$, eliminating  one color.

If in $G$,there was only one color-dominating vertex $w_s$ of color $s$ in $N(x)$, say
in $K_1$, we recolor $w_s$ by $i$. We eliminate the color $s$ by coloring each vertex of $C_s$ by a missing color.
If there is a color $t$ no more dominating we choose one, then
the color-dominating vertex $w_t$ was necessarily in $K_2$.
We color $w_t$ by $i$. We color any other vertex of $C_t$ by a missing color.
Necessarily all the remaining colors are dominating.
We conclude that $b(G-\{x\})\geq b(G)-2$.
\begin{flushright}
$\blacksquare$
\end{flushright}

Note that there exists a quasi-line graph $G$ such that $b(G-x)=b(G)-1$
for at least a vertex $x$.We give an example.
Let $\omega \geq 3$ be an integer, and let $p=2 \omega-1$.
Let $P= \{x_0,x_1,...,x_p,x_{p+1} ~\}$ be a path. We consider the graph $G_1$
obtained by replacing each edge  $[x_i,x_{i+1}]$ by a clique $K_i$ of order
$\omega$. The graph $G_1$ is a claw-free graph; we have $b(G_1)=p$ and
$b(G_1-\{x\})=b(G_1)-1$.

\section{Chordal graphs}

We want to show the following result.

\begin{theorem}
  Let $G$ be a chordal graph and $x$ be a fixed vertex of $G$. Then
  $b(G-\{x\})\geq b(G)-\omega_G $
  where $\omega_G $ is the clique number of $G$.
\end{theorem}

\vspace{4mm}
We will need first the next lemma, about the adjacencies in chordal graphs.

\begin{lemma}
Let $G$ be a chordal graph, and $a,x,b,$ be three consecutive vertices  of a cycle $\Gamma$ of $G$. Suppose that the vertex $x$ of $G$ has no neighbours in $\Gamma- \{a,b \}$. Then $a$ and $b $ 
are adjacent in $G$.
\end{lemma}
%\vspace{1cm}

\proof
The proof is by contradiction. We suppose $a$ and $b$ independent.
Suppose $\Gamma$ is a shortest cycle containing the path $axb$.
If the length of $\Gamma$ is at least $4$, then as $G$ is chordal,
and by minimality of  $\Gamma$,
it contains a chord incident with $x$ whose second endvertex is distinct from $a$ and $b$, a contradiction. 

\begin{flushright}
$\blacksquare$
\end{flushright}

In what follows, consider a $b$-coloring of $G$, with $b=b(G)$ colors. Let
$x\in V(G)$ be a fixed vertex and let $i$ be the color of $x$.
Let $I_a$ be the set of colors without color-dominating vertices in $G-\{x\}$ and
let $J_a$ be their set of color-dominating vertices in $G$. We remark that
$J_a \subset N(x)$ and no vertex in $J_a$ is neighbour of a vertex of color $i$
in $G-\{x\}$. \\

Before proving our main result, we introduce a necessary definition.

\begin{definition}
Let $x_i'$ be a fixed color-dominating vertex of color $i$, different from $x$.
Let $W=\{w\mid w\notin C_i, w ~\mbox{color-dominating vertex in} ~G\}$.
Let $k\leq b$ be a fixed integer.
We denote by $W_k$ the set of color-dominating vertices of color $k$.
A path $P$ of $G-\{x\}$ is said a {\it pseudo-alternating path} of $G-\{x\}$, and denoted by $P_k[x'_{i},z]$, if it is
a path of endvertices $x'_{i}$ and $z$, such that:

\begin{itemize}
 \item $V(P_k)\subset C_i\cup C_k\cup W \cup \{z\}$
 \item each $w\in W \cap V(P_k)$ , $w\neq x'_{i}$ and $w\neq z$, is preceded
  by a vertex of color $C_k$
 (resp. of $C_i$) and succeded by a vertex of $C_i$ (resp. of $C_k$).
 \item $V(P_k - \{z\})\cap N(x)= \varnothing$.
\end{itemize}
\end{definition}

A pseudo-alternating path $P[x_{i}',z]$ is an {\it alternating path} if $z$ is neighbour of $x$ in $G$.
We remark that if $z\notin C_i\cup C_k$, $z$ is necessarily preceded by a vertex of $C_i\cup C_k$ and, if\
$P$ is maximal, $z$ has neighbours in $C_i$ and $C_k$, belonging to $P$ or
$z\in I_a.$

%\begin{theorem}
%  Let $G$ be a chordal graph and $x$ a fixed vertex of $V(G)$. Then
%  $b(G-\{x\})\geq b(G)-\omega_G $
%\end{theorem}

{\bf Proof of the Theorem 5}

We consider a $b$-coloring of $G$, with $b(G)$ colors. Suppose $x \in C_i$.
We may assume that:\\
There is no $b$-coloring of $G-\{x\}$  by $b(G)-1$ colors \hspace{4mm} \bf(a)\\
\rm
otherwise we have the inequality of the theorem.
\rm
If $C_i-\{x\}$ contains no color-dominating vertex, we recolor each vertex of
that set  by a missing color in its neighborhood. We get a $b$-coloring of
$G-\{x\}$  by $b(G)-1$ colors, a contradiction with assumption (a).

 We may suppose from now that $C_i-\{x\}$ has color-dominating vertices.
 $x_{i}',x_{i}", \ldots, x_{i}^{(r)}$.
 Let us take $k$ in $I_a$. We recolor each vertex of $C_k$ by a missing color.
 We eliminate color $k$. In this new
 coloring, if there is a color $k'$ with no color-dominating vertex, then $k' \in I_a$.
 We recolor $C_{k'}$ and we
eliminate color $k'$. Repeating this process, we get finally a $b$-coloring by
at least $b(G)-|I_a|$.
In view to establish the bound of the theorem, we want to bound  $|I_a|$.
The bound will be established by three claims.

\vspace{2.0mm}

We denote by $Q[y,v]$ any path of $G-\{x\}$ such that $V(Q)\cap N(x)=\{v\}$.
Let $v \bar{~}$ be the neighbour of $v$ in that path.
Let $x_i$ be a color-dominating vertex.

Let $F(x_{i})$ be the set of neighbours $z$ of $x$ such that $z$ is extremity
of an alternating path $P_k[x_i, z]$ and $k$ is the color of $z$ i.e.,
$F(x_{i})~= F_1(x_{i}) \cup F_2(x_{i})$, where $F_1(x_{i})$
and $F_2(x_{i})$ are defined below:\\
\noindent
$F_1(x_{i})=\{z \in N(x)| z = z_k ~\mbox{ and in some} P_k[x_{i},z_k],
{z_k}\bar{~}~\mbox{in} C_i ~ \}$ \\
and\\
$F_2(x_{i})=\{z \in N(x)| z = z_k ~\mbox{and in some} P_k[x_{i}',z_k],
~{z_k}\bar{~} ~\mbox{in~ W} \} \diagdown F_1(x_{i}')$.\\

Let $\cal{G}$ be a component of $G-(\{x\}\cup N(x))$
Consider $X_i$ the set of color-dominating vertices of color $i$ contained in
$\cal{G}$.\\
Let $F'(x_i)= \{v \in N(x)|, \mbox {there exists~}  Q[x_{i},v] \}$.\\
Let $F'= \{v \in N(x)|, \mbox {there exists~} x_i~\rm {in}~X_i ~\rm{and}~  Q[x_{i},v]\}$.
Note that for any $x_i$, we have $F'(x_i)= F'= N(x) \cap N(\cal{G} ).$
\vspace{2mm}
Let $F_1= \cup \{F_1(x_{i}^{r}), x_{i}^{r} \in \cal{G} \}$.It is a subset of
$F'.$

\vspace{2mm}
By assumption (a) there exists a component $\cal{G}$
for which there is no recoloring of $\cal{G}$  by $\{1,...,b_G\}-\{i\}$ such
that all the colors of $W \cap({\cal{G}}~\cup N(x))$ have a color-dominating
vertex in $G-x.$ From now we use such a component.

We get the following assertion as corollary of Lemma 1.\\

\noindent \bf{Claim 1:}
%\begin{lemma}\label{lemma2}
\rm
For vertex  $x_{i}$ of  $X_i$,
 $F'(x_{i}')$ is a clique containing $F_1$.\\
 %$F'=\cup \{F'(x_{i}^{r}), x_{i}^{r} \in \cal{G} \}$ is a clique.

%\end{itemize}

%\end{lemma}

\noindent
Proof of claim 1:
%\proof
It is sufficient to note that $F'(x_i)$ is not empty, otherwise there is no
alternating path, we choose $k\neq i$ and we exchange colors $k$ and $i$ in
the pseudo-alternating paths.
No color-dominating vertex loses a color.
$X_i=\{x_{i}',\ldots, x_{i}^{r}\}$ is recolored by $k$.
A contradiction with the definition of $\cal G$.\\
$F'(x_i)$ is a clique by Lemma 1 $\bullet$
\\

At this moment we need to introduce another definition\\
Let $P_s[x_{i}',z_1]$ be an alternating path for $s\neq i$, $s\neq 1$, with
$z_1\in F_1(x_{i}')$. An {\it extension} of $P_s$, denoted by $R_s[x_{i}',z']$,
is a path of the form $P_s[x_{i}',z_1] \cup [z_1,y_i ] \cup L(y_i,z']$, where

\begin{itemize}
 \item $V(L)\subset C_i\cup C_s\cup W $; $(V(L)-\{z'\})\cap N(x)\subset W$.
 \item For each $z\in V(L)-\{z'\}$ ,
% \begin{itemize}
 \item if $z=z_k \in N(x)$, then $z \in W $ a color-dominating vertex of color $k$ and
$W_k \subset N(x)$; $z_k$ is preceded in $R_s$ by a vertex of $C_s$,
followed by a vertex $y_{i}(k)$.
 If $z\in F_{1}(x_{i}^{l})$, then $y_{i}(k)\in P_{k}(x_{i}^{l},z_k)$
\item if $z\in W - N(x)$ then $z$ is preceded by a vertex of $C_r$, followed
 by a vertex of $C_r'$, where $\{r,r'\}=\{i,s\}$
%\end{itemize}
\item  If $z'\in W$, $z'$ is preceded by a vertex of $C_{i} \cup C_{s}$.
   \\
\end{itemize}

And now, we have two claims:\\
Let $\cal K ~(A)$ be the set of colors which appear in the subset $A$ of $V$.
\noindent

\bf{Claim 2:}
\rm
\vspace{1.0mm}
%\begin{lemma}
Let $P_s[x_{i}',z']$ be an alternating extended path. Then $V(L) \cap N(x)$ is a subset of $F_1$.
So if $z'\in W$ then $z' \notin J_a.$
If $s \notin \cal{K} (F')$, then $z'\in W~-~J_a$.

%\end{lemma}

\vspace{2.0mm}
\noindent Proof of claim 2:

%\proof
Let $z_1,z_2,...,z_t,..,z_p $ be the successive vertices of $V(R_s)\cap N(x)$
and $z'=z_p.$
We show by induction on $t$, that $z_t$ is in $F_1$.
Suppose that $z_{t-1}\in F_{1}(x_{i}^{l})$ for some $l$. So there is a path
$P_{t-1}[x_{i}^{l},y_{i}^{(t-1)}]$.
Composing it with $R_{s}(y_{r}^{(t-1)}, z_t)$, we get a path $Q[x^{(l)}, z_t]$.
If $z_t \notin F_1$, we do $\tau (i,t)$ from $x_{i}^{l'}$ for any $l'$, where
$\tau (i,t)$ means exchanging the colors $i$ and $t$ in $\cal G$. No color-dominating vertex
loses color $t$ even if it is an  extremity of an alternating path $P_t$ in
this later case it is neighbour of $z_t$.
Some color-dominating neighbours of
$N_t(x)$ may lose color $i$. We recolor each $v\in X'_i$ by a missing color.
No color-dominating vertex of color $i$ is created by uniqueness of the color-dominating of color $t$.
We get a coloring by $b(G)-1$ colors. A contradiction. So $z_t \in F_1$ and $z_t \notin J_1$.
If $s \notin \cal{K} (F')$, then no vertex of $C_s$ is in $F_1$. So $z'\in W$.
As $F_1 \cap J_a=\emptyset$, then $z'\in W-J_a .$
%\begin{flushright}
%$\blacksquare$
%\end{flushright}
\noindent

\vspace{2.0mm}
\bf{Claim 3:}
\rm
%\begin{lemma}\label{F'}
$\cal{K}$$(F')$ contains $I_a$

%\end{lemma}
\vspace{2.0mm}

\noindent Proof of claim 3:
\vspace{2.0mm}

It is by contradiction. We suppose that there is a color $s$ such that
$s \in I_a \setminus \cal{K} (F')$. Thus no vertex of color $s$ belongs to $F_1$.

\vspace{2mm}
{\it Case 1:} There is no path $P_s[x_{i}^{r},z]$ with $x_{i}^{r}\in \cal{G}$ and
$z\in N(x)$.\\
We do $\tau (i,s)$ along the pseudo-alternating paths $P_s[x_{i}^{r},u]$,
$u\notin N(x)$ simultaneously.
No color-dominating vertex loses color.
We recolor the remaining vertices of $C_i$ in $\cal{G}$ and this leads to a
contradiction with the assumption.

\vspace{2mm}
{\it Case 2:} There exists $x_{i}^{r}$ in $\cal{G}$ and a path $P_s[x_{i}^{r},z_t]$ with $z_t \in N(x)$, $x_{i}^{r} \in \cal{G}$.

\vspace{2mm}
This case will be divided in two sub-cases:

\vspace{2mm}
Case 2.1: There exists $ P_s[x_{i}^{r},z_t]$ with $z_t \notin F_1$

\vspace{2mm}
By claim 1, $F'$ does not contain $z'_{t}$, with $z'_{t}\neq z_{t}$. We do $\tau (i,t)$ simultaneously in all alternating $P_t[x_{i}^{l},z]$, with $x_{i}^{l}$ in $\cal{G}$.
As $z_t \in F'$,by Claim1,
the color-dominating vertices contained in $F'$ do not lose color $t$,
they may lose color $i$. The color-dominating vertices which may lose a color are the
color-dominating vertices preterminal in $P_t[x_{i}^{l},z_t]$; these color-dominating vertices may lose
color $i$. We recolor $C_i$ by missing colors in $\cal{G}$.

\vspace{2mm}
Case 2.2: For any $P_s$, the extremity contained in $N(x)$ is in $F_1$.\\
So this extremity does not belong to $J_a$.\\
If for any $x_{i}^{l} \in X_i$, there is no alternating path $P_s[x_{i}^{l},z]$
with color $s$ preceding $z$, then we do
$\tau (i,s)$ in all the pseudo-alternating paths and the paths $P_s$.
We have the same conclusion as in case 2.1.

\vspace{2mm}
If for some $l, z$ is preceded by a vertex of color $s$ in $P_s[x_{i}^{l},z]$
we consider the extended paths $R_s$. We do
$\tau (i,s)$ simultaneously in all alternating and pseudo-alternating paths
$P_s$ and the extended paths $R_s$.
The color-dominating vertices which may lose a color are either in $N(x)$ or
in $N^{2}(x)$, they are among terminal vertices and preterminal vertices of
the alternating paths $P_s$ and $R_s$; and they may lose color $i$. We then
recolor each remaining vertex of $C_i \cap {\cal G}$ by a missing color.

\vspace{2.0mm}
In each case we have a contradiction with the definition of the component
$\cal G$. So $I_a \subset \cal K$$(F') \bullet$ 
\vspace{2.0mm}

By claim 1 and claim 3 there is a clique containing $x$ and $F'$. So we have
$\omega(G)\geq \vert F'\vert+1  \geq \vert I_a \vert +1$, and this finishes the
proof of the theorem.
\begin{flushright}
$\blacksquare$
\end{flushright}

\section{Graphs with large girth}

%Let $x_1,...,x_p$ be a sequence of vertices of a graph $G$ such that\\
%$d(x_1) \leq ... \leq d(x_i)\leq d(x_{i+1})\leq...\leq d(x_{p}).$

The \textit{m-degree} of a graph $G$, denoted by $m(G)$, is the largest integer $m$ such that $G$ has $m$ vertices of degree at least $m-1$. It is known that, for any graph $G$, $b(G) \leq m(G)$ (see  \cite{Irving}).
%Let $m(G)=max \{p, d(x_i)\geq (i-1)  \mbox{for each}  i\leq p \}.$

Note that A.Campos et al. \cite{Camp} have shown that graphs of girth at least
$7$ have high b-chromatic number; for each such a graph $G$ this number is at
least $m(G)-1$. 

We can verify that $m(G-\{x\}) \geq m(G)-1.$ Indeed, we have three possibilities to consider:
$x$ is one of the $m$ vertices of degree $m-1$; $x$ is a neighbor of one of the $m$ vertices of degree $m-1$, or $x$ is not in any of the previous situations. In the first case, there remain $m-1$ vertices of degree at least $m-2$ and, thus $m(G-\{x\}) \geq m(G)-1.$ In the second case, there remain $m$ vertices of degree at least $m-2$, and again, $m(G-\{x\})\geq  m(G)-1.$ In the latter case the m-degree does not change, that is, $m(G-\{x\})= m(G)$.

So $b(G-\{x\}) \geq b(G)-2$ for graphs of girth at least $7$. No particular bound is known for graphs of girth $5$ or $6$.
In this work we show the following.

\begin{theorem}
 Let $G$ be a graph of girth al least $5$.
 For each vertex $x$, \\
 $b(G-\{x\}) \geq b(G)-2$
\end{theorem}
\proof
Let $\cal B$ be a b-coloring of $G$ and let $i$ be the color of the deleted
vertex $x$.
Let $W$ be the set of color-dominating vertices of colors different from $i$ in $G$.
Let $W_k$ be the subset of those of color $k$.\\
We may suppose that there is a set of  color-dominating vertices
$X_i$ of color $i$, different from $x$.
For each vertex $u$ let $K_1(u)$ be the set of colors with at least a color-dominating vertex
in $N(u)$, let us set $K_2(u)= \{1,2,...,b\}- (K_1(u) \cup \lbrace i\rbrace).$

We use the notations of the previous section. We may suppose $|I_1| \geq 3.$
Note that for $u\neq x$, as the girth is at least $5$, $K_1(u)$ does not
contain $I_1$. So $K_2(u)$ is not empty as it intersects $I_1$.
By definition of $I_1$, the color $i$ is a missing color for
each vertex of $J_1$ in $G- \{x\}.$ 
So for any $x' \in X_i$, $K_1(x')$ does not intersect $I_1$.

(a)  Let $x' \in X_i$. Let $k \in K_2(x')$ be fixed.
\\
($\bf P_a$) \rm
If $G-N(N_k(x'))$ intersects $W_p$ for each color $p$ different from $i$,
\rm
\\
we color $x'$ by $k$, each vertex of $N_k(x')$ by a missing color.
So $|X_i|$ decreases.
The color-dominating vertices of $K_1(x')$ may lose color $i$ in their neighborhood.\\
\indent
(b) As long as there exists a vertex $x''$ of $X_i$ satisfying ($\bf P_a$)
\rm
we do a recoloring.
\\
From now we may suppose that $X_i$ is not empty and no vertex of $X_i$
satisfies (a).

\begin{lemma}
  Let $x'$ be a fixed vertex of $X_i$.
  %Let $K_1(x')$ be the set of colors with at least a color-dominating vertex in $N(x')$.
  %Let $N'$ be the subset $N(x_i)$ not colored by $K_1$.
  For each $k \in K_2(x')$ there exists exactly one $j_k$ such that $W_{j_k}$
  is contained in $N(N_k(x'))$ and $N(N_k(x')) \cap W_j=\emptyset$ for
  any other $j$.

\end{lemma}
\proof
We know that the property $(a)$ is not satisfied. As $g(G)\geq5$, if a vertex
$w_j$ is in $N(N_k(x'))$ then $w_j$ is not neighbour of $N_r(x')$ for $r$
different from $k$.
As $(a)$ is not satisfied, it follows that for each $k\in K_2(x')$, $N_k(x')$ is
neighbour of any vertex of a set $W_j$ for some $j$. As $g(G) \geq 5$, $j$ is in
$K_2(x')$ and $W_j$ is not neighbour of $N_t(x')$ for $t \neq k$.
It follows that  for $k$ fixed, $j$ is unique.
This finishes the proof of the lemma.
\\
\vspace{2mm}

Let $s $ be a fixed element of $I_1$. Let $x' \in X_i$.
We know that the set $I_1$, by definition, is a subset of $K_2(x').$ By the
precedent Lemma there is exactly one $k$ such that $W_s \subset N(N_k(x'))$.
We color each vertex of $N_k(x')$ by a missing color
and $x'$ by $k$. If $N_k(x')$ meets $N_k(x'')$ for some $x'' \in X_i$, by the
precedent Lemma, we have  $W_s \subset N(N_k(x''))$; we color $x''$ by $k$ as
well and we recolor $ N_k(x'')$ by missing colors different from $i$.
We recolor so each vertex of $X_i$. Then we recolor each vertex of color $i$
by a missing color.
If,finally,the color $s$ has no vertex dominating the colors $\{1,..,b\}-\{i\}$,
we recolor each vertex $u_s$ by a missing color different from $i$.
\\

After this recoloring of color $i$ and eventually color $s$, we get a b-coloring of $G$ by at least $b(G)-2$ colors.
\begin{flushright}
$\blacksquare$
\end{flushright}

\vspace{5mm}

The first author acknowledges partial support by CAPES and CNPq.

\vspace{7.0mm}

\noindent Renata Del-Vecchio\\
renata@vm.uff.br\\
Instituto de Matematica,\\ 
Universidade Federal Fluminense\\
Niteroi, RJ, Brazil\\
\vspace{8.0mm}

\noindent Mekkia Kouider\\
km@lri.fr\\
Universite Paris-Sud\\
Paris, France.

\end{document}